\documentstyle[12pt]{amsart}

\newtheorem{definition}{Definition}
\def\bbeta{\gamma}
\def\M{$\Sigma$}
\def\S{\Sigma}

\title{ Noether's  variational theorem II and the BV formalism}
\begin{document}

\author{Ron Fulp}
\address{Department of Mathematics, North Carolina State University,
Raleigh NC 27695}
\email{fulp\@math.ncsu.edu}
\author{Tom Lada}
\address{Department of Mathematics, North Carolina State University,
Raleigh NC 27695}
\email{lada\@math.ncsu.edu}
\author{Jim Stasheff }
\address{Department of Mathematics, University of North Carolina, Chapel
Hill, NC 27599-3250, USA}
\email{jds\@math.unc.edu}
\thanks{Stasheff's research supported in part by the NSF throughout most
of his
career, most recently under grant DMS-9803435. }

\smallskip

\maketitle

\section{Introduction}

Lagrangian physics derives `equations of motion' from a variational
principle of least action.  Here an action refers to an integral
$$
S(\phi)=\int_M L((j^n\phi)(x)) vol_M
$$
over some manifold $M$ where $\phi$ is a (possibly vector valued)
function
on $M$
or section of a bundle $E$ over $M$ and $L$ is a `local function' on
$E$,
meaning a
function on some finite jet space $J^n E.$

The Euler-Lagrange equations describe the
critical points of $S$ with respect to variations in $\phi.$ The action
may
have symmetries, i.e. variations in $\phi$ which do not change the value
of
$S$ and hence are physically irrelevant in the sense that $\phi$ and its

transformed
value encode the same physical  information.
Noether's second variational theorem establishes a correspondence
between
symmetries and
differential algebraic relations among the Euler-Lagrange equations.

These symmetries create difficulties for quantization of such physical
theories.  The method
of Batalin and Vilkovisky \cite{bv:closure,bv:closure2} was invented
to handle these difficulties, but
turns out to also be of
interest in a classical context.  Their method extends the BRST
cohomological approach by
introducing anti-fields (independently and previously due to Zinn-Justin
\cite{zj,zj:leshouches})  dual to the original fields and anti-ghosts
which
correspond to the
Noether relations and are dual to the ghosts which generate the BRST
complex. A key
ingredient in their approach is to use the duality to give an
anti- bracket (independently due to Zinn-Justin \cite{zj,zj:leshouches}
and also
known as an odd Poisson or Gerstenhaber bracket) in their construction.

The relevance of Noether's theorem is not emphasized in most of the
literature using the BV approach. One aim of the present paper is to
restore
such an emphasis: Part of the BV complex is the Koszul-Tate resolution
of the differential ideal generated by Euler-Lagrange equations. The
anti-fields
generate the Koszul complex, which is not a resolution; the anti-ghosts
provide
the next level of generators as described by Tate corresponding to the
relations among the
Euler-Lagrange equations. Rather than carrying out this analysis in the
abstract, we illustrate it explicitly in terms of the Poisson sigma
models
of Cattaneo and Felder.

The higher order terms in the BV differential in these examples
can be related to parts of an $L_\infty$-algebra structure, as we will
explain elsewhere.

In Section 2, we review the basics of the Lagrangian approach and
establish
the notation we will use.  Section 3 is devoted to Noether's Second
Theorem
with a slight
modernization of language and notation. In Section 4, we present the
Cattaneo-Felder
sigma model and work out the Noether identities. In Section 5, we begin
the
description
of the Batalin-Vilkovisky formalism, pointing out the initial
Chevalley-Eilenberg (or
BRST) part of the differential and especially the Koszul-Tate part. The
latter shows
explicitly how the anti-ghosts encode the Noether identities. We also
recall how to
extend the gauge symmetries to act on the anti-fields and anti-ghosts.
To
combine the
Koszul-Tate and Chevalley-Eilenberg differentials into a total
differential
of square
zero requires `terms of higher order', which are created via the
Batalin-Vilkovisky
anti-bracket as worked out in  Section 6.

\section{Preliminaries}

Let $\S$ be an $s$-dimensional manifold and $\pi:E\rightarrow \Sigma$ a
vector
bundle of fiber dimension $k$ over $\S$. Let $J^{\infty}E$ denote the
infinite jet bundle of $E$ over $\S$ with $\pi^{\infty}_E:J^{\infty}E
\rightarrow E$ and $\pi^{\infty}_\S:J^{\infty}E \rightarrow \S$ the
canonical projections. The vector space of smooth sections of $E$ with
compact support will be denoted $\Gamma E$. For each section
$\phi$ of $E$, let $j^{\infty}\phi$ denote the induced section
of the infinite jet bundle $J^{\infty}E$. We will consider `local'
functions defined on a finite jet space (see below), but  refer to
$J^{\infty}E$ to avoid specifying some finite jet.

The restriction of the infinite jet bundle over  an appropriate open
$U\subset \Sigma$ is trivial with fibre an infinite dimensional vector
space
$V^\infty$.  The bundle
\begin{eqnarray}
\pi^\infty : J^\infty E_U=U\times V^\infty \rightarrow U
\end{eqnarray}
then has induced coordinates given by
\begin{eqnarray}
(x^i,u^a,u^a_i,u^a_{i_1i_2},\dots,).
\end{eqnarray}
We use multi-index notation and the summation convention throughout the
paper. If $j^{\infty}\phi$ is the section of $J^{\infty}E$ induced by
a section $\phi$ of the bundle $E$, then $u^a\circ
j^{\infty}\phi=u^a\circ \phi$ and $$u^a_I\circ j^{\infty}\phi=
(\partial_{i_1}\partial_{i_2}...\partial_{i_r})(u^a\circ
j^{\infty}\phi)$$ where $r$ is the order of the symmetric multi-index
$I=\{i_1,i_2,...,i_r\}$,with the convention that, for $r=0$, there are
no derivatives.

\begin{definition} We say that a real-valued function
on the jet space $J^{\infty}E$ is a {\bf local function} if it is the
composite
of the projection from $J^{\infty}E$ onto $J^kE$ and a smooth
real-valued function on  $J^kE$ for some k. Thus such functions
are pull-backs of functions in $C^{\infty}(J^kE)$ under the projection
$\pi^{\infty}_k: J^{\infty}E\longrightarrow J^kE.$
\end{definition}
\smallskip

Let
\begin{eqnarray}
D_i=\frac {\partial}{\partial x^i}+u^a_{iJ}\frac {\partial}{\partial
u^a_J}
\end{eqnarray}
be the total differential operator acting on the space $Loc_E$ of local
functions
defined on the jet space $J^{\infty}E.$

More generally, {\bf total differential operators} are mappings from
$Loc_E$ into $Loc_E$ defined in local coordinates by $Z=Z^ID_I$
where $Z^I\in Loc_E$ and
and $D_I=D_{i_1}\circ \cdots \circ D_{i_r}$ for each symmetric
multi-index $I.$

It can be shown that the complex
$\Omega^*(J^{\infty}E,d)$ of differential forms splits as a bicomplex
(though the finite
level complexes $\Omega^*(J^pE)$ do not). The bigrading is described
by writing a differential $p$-form
$\alpha =\alpha_{IA}^{\bf J}(\theta^A_{\bf J}\wedge dx^I)$ as an element
of
$\Omega^{r,t}(J^{\infty}E)$, with $p=r+t,$ where
\begin{eqnarray} dx^I=dx^{i_1}\wedge...\wedge dx^{i_r}, \quad \quad
  \quad \theta^A_{\bf J}=\theta^{a_1}_{J_1}\wedge...\wedge
  \theta^{a_t}_{J_t} \end{eqnarray}
and  $$\theta^a_J=du^a_J-u^a_{J\mu}dx^{\mu}.$$
We restrict the complex $\Omega^*$ by requiring that the
functions $\alpha_{IA}^{\bf J}$ be local functions.
In this context,
the horizontal differential is obtained by noting that $d\alpha$ is in
$\Omega^{r+1,t}\oplus \Omega^{r,t+1}$ and then denoting the two pieces
by, respectively, $d_H \alpha$ and $d_V\alpha$.

We will work exclusively with the $d_H$ subcomplex, the algebra of
horizontal forms $\Omega^{*,0}$, which is the exterior algebra in the
$dx^i$ with coefficients that are local functions. In this case we
often use Olver's notation $D$ for the horizontal differential
$d_H=dx^iD_i$ where $D_i$ is the total derivative defined above.
It is well-known that in this language, the Poincar\'e lemma asserts
that
on an appropriate open subset of $J^{\infty}E,$
$d_H\alpha=0$ for $\alpha \in \Omega^{s,0}$ iff
$$\alpha=\partial_{\mu} j^{\mu}(dx^1\wedge \cdots \wedge dx^s)$$ for
some
choice
of local functions $\{j^{\mu}\}.$

\begin{definition}
A {\bf local functional} is a function $S$ from $\Gamma E$ into the
reals such that, for each section $\phi\in \Gamma E,$ we have
\begin{eqnarray}
S(\phi)=\int_\Sigma L(x,\phi^{(p)}(x)) dvol_\Sigma = \int_\Sigma
(j^\infty
\phi)^*  L(x,u^{(p)}) dvol_\Sigma
\end{eqnarray}
is the integral over $\S$ of the pull-back $(j^{\infty}\phi^*)L$of some
local
function $L$ on $J^{\infty}E.$ Recall that the elements of $\Gamma E$
have
compact support so that the integral is well-defined .
\end{definition}

These definitions  reflect the fact that we identify the fields
$\phi$ of a physical theory with sections of an appropriate vector
bundle $E\longrightarrow \Sigma.$ With this identification, the
{\em Lagrangian} $L$ of the theory is a local function on
$J^{\infty}E.$
We work on $J^{\infty}E$ for convenience but the Lagrangian, being
local, only depends on finitely many derivatives of the fields. Finally,

the {\em action} $S$ corresponding to  the Lagrangian is simply the
local
functional defined by $L$ as in the definition above.

\begin{definition}
{\em The Euler-Lagrange operator}: For $1\leq a\leq k$,
let $E_a$ denote the $a$-th component of the {\em Euler-Lagrange
operator}
defined for $F\in Loc_E$ by
\begin{eqnarray}
E_a(F)=\frac{\partial F}{\partial u^a}-\partial_i\frac{\partial F}
{\partial u^a_i}+\partial_i\partial_j\frac{\partial F}{\partial
u^a_{ij}}-...=(-D)_I(\frac{\partial F}{\partial u^a_I}).
\end{eqnarray}
\end{definition}

We say that $Q$ is an {\em evolutionary vector field on E} if it is a
mapping from $J^{\infty}E$ into the vertical vector fields on $E.$
In local coordinates $Q=Q^a \frac {\partial}{\partial u^a}$
where the functions $Q^a$ are local functions on $J^{\infty}E.$
For every evolutionary vector field $Q$ on
$E$, there exists its {\em prolongation},denoted $pr(Q),$the unique
vector
field on $J^{\infty}E$ such that
$(d\pi^{\infty}_E)(pr(Q))=Q$ and ${\cal L}_{pr(Q)}({ C})\subseteq {
C}.$ Here ${\cal L}_{pr(Q)}$ denotes the Lie
derivative operator with respect to the vector field $pr(Q).$ The ideal
 ${C}$ is the ideal of forms on $J^{\infty}E$ generated by
the contact forms $\{\theta^a_J\}$ used above in the definition
of the bicomplex.

In local adapted coordinates, the prolongation of an evolutionary
vector field $Q=Q^a\partial/\partial u^a$ assumes the form
$pr(Q)=(D_JQ^a)\partial/\partial u^a_J$.

Given a total differential operator $Z,$  define a new total
differential
operator
$Z^+$ called the (formal) {\em adjoint} of $Z$ by
\begin{eqnarray}
\int_M (j^{\infty}\phi)^*(FZ(G))dvol_{\S}=\int_{\S}
(j^{\infty}\phi)^*(Z^+(F)G)dvol_{\S}
\end{eqnarray}
for all sections $\phi\in \Gamma E$ and all $F,G\in Loc_E$. It
follows that
\begin{eqnarray}
FZ(G)dvol_{\S}=Z^+(F)Gdvol_{\S} + d_H\zeta\label{x}
\end{eqnarray} for some
$\zeta\in \Omega^{n-1,0}(E)$. If $Z=Z^JD_J$ in local coordinates, then
$Z^+(F)=(-D)_J(Z^JF)$. This follows from an integration by parts in
(\ref{x}) and the fact that (\ref{x}) must hold for all $G.$

\section{Gauge symmetries and Noether identities}

Recall that if a Lie group G acts as automorphisms of a vector bundle
$E\longrightarrow \Sigma$
(over the identity of \M) in such a way that it leaves the action S
of a Lagrangian $L:J^{\infty}E\longrightarrow {\bf R}$  invariant, then
the group action induces a vertical vector field  $\tilde \eta$ on E,
for each
element $\eta$ of the Lie algebra of G, such that the prolongation
$pr(\tilde \eta)$ of $\tilde \eta$ to $J^{\infty}E$ has the property
that
$dL(pr(\tilde \eta)) vol_M$ is $d_H$ exact.
Here $vol_\Sigma$ denotes both a
volume on $\S$ and its pullback to $J^{\infty}E$ via
the projection $J^{\infty}E\longrightarrow \Sigma.$ An evolutionary
vector
field
$Q_E$ on E is called a {\em variational symmetry}
of a Lagrangian L iff it has the property that
$dL(pr(Q_E)) vol_\Sigma=pr(Q_E)(L) vol_\S$ is
$d_H$ exact. In local coordinates, $Q_E$ is a variational symmetry iff
 \begin{equation}
pr(Q_E)(L):=D_K(Q_E^a)\frac  {\partial L}{\partial u^a_K}
\end{equation}
is a divergence, i.e., iff it is equal to $D_{\mu}j^{\mu}$ for some set
$\{j^{\mu}\}$ of local functions defined on  $J^{\infty}E.$
``Integrating
by parts" shows that this condition is equivalent to requiring that
$Q_E^a(-D)_K(\frac  {\partial L}{\partial u^a_K}))$ be a divergence.
But the Euler Lagrange operator $E_a$ acting on the Lagrangian $L$ is
defined
by the equation $E_a(L)=(-D)_K(\frac  {\partial L}{\partial u^a_K}).$
Thus an evolutionary vector field $Q_E$ is a variational symmetry of a
Lagrangian $L$ iff $Q_E^aE_a(L)$ is a divergence.

Finally, a gauge
symmetry of a Lagrangian $L$ is defined when there is a linear mapping
from
$Loc_E$
into the
variational symmetries. To be more precise, there must exist local
functions
$R^{aI}:J^{\infty}E\longrightarrow {\bf R}$ such that
$R^{aI}(D_I\epsilon)\frac {\partial}{\partial u^a}$ is a variational
symmetry
of $L$ for each local function $\epsilon:J^{\infty}E\longrightarrow {\bf
R}.$
Notice that the coefficients of the vector field depend linearly on both
$\epsilon$ and its derivatives. It follows that being a gauge symmetry
is equivalent to requiring that $R^{aI}(D_I\epsilon)E_a(L)$ be a divergence
for each $\epsilon.$ This in turn is equivalent to saying that
$\epsilon (R^{aI}D_I)^+(E_a(L))$ is a divergence for each $\epsilon.$
Here
$(R^{aI}D_I)^+$ is the formal adjoint of the differential operator
$R^{aI}D_I$
which was defined in Section 2. The adjoint of a differential operator
is
also a differential operator and consequently there exist local
functions
$R^{+aI}:J^{\infty}E\longrightarrow {\bf R}$ such that
$(R^{aI}D_I)^+=R^{+aI}D_I.$
These functions are found by working out the iterated total derivatives
$(-D)_I(R^{aI}F).$

In many cases it is easier to use an ``integration by parts"
procedure to obtain the coefficients $\{R^{+aI}\}.$ This is what we do
for the
Poisson $\sigma$-model below.

It follows easily that
$\epsilon \mapsto R^{aI}(D_I\epsilon)\frac {\partial}{\partial u^a}$
defines a gauge symmetry iff $\epsilon R^{+aI}D_I(E_a(L))$ is a
divergence for
each $\epsilon.$ Finally, this condition is equivalent to saying that
$R^{+aI}D_I(E_a(L))$ is identically zero on the jet bundle. Such
identities
are called {\em Noether identities} or {\em dependencies} in the
translation of
Noether's original term. One thus has a one-one correspondence between
gauge symmetries of a Lagrangian and Noether identities.

The original version of Noether
in `Invariant variation problems' \cite{noether}, was written in terms
of an
infinite continuous group, $G_{\infty \rho},$ `understood to be a group
whose most
general transformations depend on $\rho$ essential arbitrary functions
and
their
derivatives'. Noether's Theorem II refers to an integral $I$ ($= S$ in
our
notation) and
reads:
\begin{quote}If the integral I is invariant with respect to a $G_{\infty
\rho}$
in which the arbitrary functions occur up to the $\sigma$-th derivative,
there
there subsist $\rho$ identity relationships between the Lagrange
expressions and
their derivatives up to the $\sigma$-th order. $\dots$ the converse
holds.
\end{quote}
Later in that paper these relations are called {\em dependencies}.
\vskip2ex
To recast and summarize in our notation and terminology, we have:
\vskip2ex
\noindent THEOREM (Noether) For a given Lagrangian $L$ defined on the
jet
bundle
$J^{\infty}E$ and for local real-valued functions $\{R^{aI}\}$ defined
on
$J^{\infty}E,$ the following statements are equivalent:\newline
(1) The functions $\{R^{aI}\}$ define a gauge symmetry of
$L,$ i.e., $R^{aI}(D_I\epsilon)\frac {\partial}{\partial u^a}$ is a
variational
symmetry of $L$ for each local function
$\epsilon:J^{\infty}E\longrightarrow {\bf R}.$
\newline
(2) $R^{aI}(D_I\epsilon)E_a(L)$ is a divergence for each
$\epsilon.$\newline
(3) The functions $\{R^{aI}\}$ define Noether identities of
$L,$ i.e.,  $R^{+aI}D_I(E_a(L))$ is identically zero on the jet bundle.

\section{The Poisson sigma model}
\def\aaa{\alpha}

To provide a specific example of this correspondence and how
it relates to the Batalin-Vilkovisky machinery, we turn to a
Poisson sigma model of Cattaneo and Felder \cite{cf:defquant}.

The fields of this Poisson $\sigma$-model are ordered pairs $(X,\eta)$
such that $X$ is a mapping from a 2-dimensional manifold $\Sigma$ into a

Poisson
manifold $M$ and $\eta$ is a section of the bundle
$Hom(T\Sigma,X^*T^*M)\longrightarrow \Sigma.$
These fields are subject
to boundary conditions, namely they should satisfy the conditions:
$X(u)=0$ and
$\eta(u)(v)=0$ for arbitrary $u$ in the boundary of $\Sigma$ and for $v$ tangent to the
boundary of $\Sigma$ at $u.$  Observe that for
each $u\in \Sigma,$ we can regard
$\eta(u)$ as a linear mapping from $T_u\Sigma$ into $T^*_{X(u)}M.$ In
local
coordinates $\{u^{\mu}\}$ on $\Sigma$ and $\{x^i\}$ on $M,$ we write
$dX=(dX^j) \frac {\partial}{\partial x^j}$ and
$\eta(\frac {\partial}{\partial u^{\mu}})=\eta_{i,\mu}dx^i.$
The Poisson structure is given
by a {\em Poisson tensor} which is a skew-symmetric tensor on $M$
\begin{equation}
 \alpha = \aaa^{ij}(\frac {\partial}{\partial x^i} \wedge
\frac{\partial}{\partial x^j})
\end{equation}
which satisfies a Jacobi condition:
\begin{equation}\label{e-Jacobi}
\alpha^{il}\partial_l\alpha^{jk}+
\alpha^{jl}\partial_l\alpha^{ki}+
\alpha^{kl}\partial_l\alpha^{ij}=0,
\end{equation}
The action $S$ of the model is defined in such local coordinates by
\begin{equation}
S(X,\eta)=\int_{\Sigma}(\eta_i\wedge dX^i) +
\frac {1}{2} (\alpha^{ij}\circ X)(\eta_i\wedge \eta_j).
\end{equation}
To understand this action in a more invariant notation, recall that for
each $u\in
\Sigma,$ $dX$ is a linear mapping from $T_u\Sigma$ into $T_{X(u)}M$ and
so
one may define a two-form $\eta\wedge dX$ on $\Sigma$ by
\begin{equation}
(\eta\wedge dX)(v_1,v_2)=\eta(v_1)(dX(v_2))-\eta(v_2)(dX(v_1))
\end{equation}
for $v_1,v_2\in T\Sigma.$ We may also define a two-form
$\alpha_X(\eta\wedge\eta)$ on $\Sigma$ by
\begin{equation}
\alpha_X(\eta\wedge\eta)(v_1,v_2)=(\alpha\circ X)(\eta(v_1),\eta(v_2)).
\end{equation}
Using the coordinates defined above, we see that:
\begin{eqnarray}
\eta\wedge dX=\eta_i\wedge dX^i \quad \quad \\
\alpha_X(\eta\wedge\eta)=(\alpha\circ X)(\eta\wedge\eta) =
\frac {1}{2} \alpha_X^{ij}(\eta_i\wedge \eta_j).
\end{eqnarray}

{\em For the remainder of the paper, we will restrict to $M={\mathbf
R}^k$
to avoid inserting `in local coordinates' repeatedly.}

According to the variational principle, we obtain extrema of $S$ as
those
fields
$(X,\eta)$ which satisfy the Euler-Lagrange equations:
\begin{equation}
E_{X^i} := d\eta_i+\frac{1}{2}\partial_i\alpha^{jk}(\eta_j\wedge
\eta_k)=0
\label{Eulerform1}\end{equation}
and
\begin{equation}
E_{\eta_i} :=-dX^i-\alpha^{ij}\eta_j=0.
\label{Eulerform2}\end{equation}

In terms of the components of the fields, we write
\begin{equation}
E_{X^i}=(\partial_{\mu}\eta_{i,\nu}+
\frac{1}{2}\partial_i\alpha^{jk}\eta_{j,\mu}\eta_{k,\nu})\epsilon^{\mu\nu}
\end{equation}
and
\begin{equation}
E_{\eta_{i,\nu}}=-(\partial_{\mu}X^i+\alpha^{ij}\eta_{j,\mu})\epsilon^{\mu\nu}.
\end{equation}
The gauge symmetries of the action are parameterized by  all sections
$\beta$ of the
bundle $X^*T^*M\longrightarrow \Sigma$ which vanish on the boundary of
$\Sigma.$
 For each such $\beta,$ define
$\delta_{\beta}$ acting on the fields by
\begin{equation}
(\delta_{\beta}X)^i=(\alpha\circ X)(dx^i,\beta)
\end{equation}
\begin{equation}
(\delta_{\beta}\eta)(W\circ X)=-(d\beta)(W\circ X)-
(({\cal L}_W\alpha)\circ X)(\eta,\beta)
\end{equation}
where   $W$ is a vector field on $M,$ and  ${\cal
L}_W\alpha$ is the Lie derivative of $\alpha$ with respect to $W.$
Observe that $\delta_{\beta}X$ and $\delta_{\beta}\eta$ are indeed
again fields since $\delta_{\beta}X$ is a mapping from $\Sigma$ to $M$
and $\delta_{\beta}\eta$ is a section of the bundle
$Hom(T\Sigma,X^*T^*M)\longrightarrow \Sigma.$

If we regard
$X^i$ and
$\eta_{i,\nu}$ as jet coordinates on an appropriate
jet bundle, we may write $\delta_{\beta}$ as a variational symmetry
\begin{equation}
\delta_{\beta}=(\alpha_X^{ij}\beta_j)\frac {\partial}{\partial X^i}
-(\partial_{\mu}\beta_i+((\partial_i\alpha^{jk})\circ
X)\eta_{j,\mu}\beta_k)\frac {\partial}{\partial \eta_{i,\mu}}.
\end{equation}
For notational convenience, we will not show the explicit
$X$ dependence throughout the remainder of this section except
when it is misleading to fail to do so.

It follows from Noether's theorem that
\begin{equation}
(\alpha^{ij}\beta_j)E_{X^i}
-(\partial_{\mu}\beta_i+\partial_i\alpha^{jk}\eta_{j,\mu}\beta_k)E_{\eta_{i,\mu}
}
\end{equation}
is a divergence. To find the corresponding Noether identity, we must be
able to
factor out the gauge parameters $\beta_k,$ so we transform the term
$(\partial_{\mu}\beta_i)E_{\eta_{i,\mu}}$  via the identity
\begin{equation}
(\partial_{\mu}\beta_i)E_{\eta_{i,\mu}}=\partial_{\mu}(\beta_iE_{\eta_{i,\mu}})
-\beta_i\partial_{\mu}E_{\eta_{i,\mu}}=div
-\beta_i\partial_{\mu}E_{\eta_{i,\mu}}.
\end{equation}
If
$N^k :=\alpha^{ik}E_{X^i}
+\partial_{\mu}E_{\eta_{k,\mu}}-
\partial_i\alpha^{jk}\eta_{j,\mu}E_{\eta_{i,\mu}},$ then
\begin{equation}
N^k\beta_k=(\alpha^{ik}E_{X^i}
+\partial_{\mu}E_{\eta_{k,\mu}}-
\partial_i\alpha^{jk}\eta_{j,\mu}E_{\eta_{i,\mu}})\beta_k
\end{equation}
is a divergence for every $\beta_k.$ It follows that the integral of
$N^k\beta_k$ vanishes for all $\beta_k$ and consequently $N^k=0$ for
each
k. From this, we see that the equations
\begin{equation}
\alpha^{ik}E_{X^i} +\partial_{\mu}E_{\eta_{k,\mu}}-
\partial_i\alpha^{jk}\eta_{j,\mu}E_{\eta_{i,\mu}}=0
\end{equation}
are the Noether identities corresponding to the gauge symmetry
$\delta_{\beta}$
defined above.

To write this identity in differential form notation, multiply the last
equation
by $du^1\wedge du^2$ and use the identity $\epsilon^{\mu\nu}(du^1\wedge
du^2)=du^{\mu}\wedge du^{\nu}$  to get
$$
\alpha^{ij}(\partial_{\mu}\eta_{i,\nu}
+\frac{1}{2}\partial_i\alpha^{rs}\eta_{r,\mu}
\eta_{s,\nu})(du^{\mu}\wedge
du^{\nu})
-\partial_i\alpha^{kj}\eta_{k,\mu}(\partial_{\nu}X^i
+\alpha^{ir}\eta_{r,\nu})(du^{\mu}\wedge du^{\nu})
$$
$$
+\partial_{\mu}(\partial_{\nu}X^j
+\alpha^{jr}\eta_{r,\nu})(du^{\mu}\wedge du^{\nu})=0,$$
which in turn implies that
\begin{equation}
\alpha^{ij}[d\eta_i+\frac{1}{2}\partial_i\alpha^{rs}(\eta_r\wedge
\eta_s)]
-\partial_i\alpha^{kj}[\eta_k\wedge
(dX^i+\alpha^{ir}\eta_r)]+d[dX^j+\alpha^{ji}\eta_i]=0.
\end{equation}

Now utilize the formulas  (\ref{Eulerform1})
 and (\ref{Eulerform2})
for $E_{X^i}$ and $E_{\eta_{i,\mu}}$ above to obtain the Noether
identities in the form
\begin{equation}
\alpha^{ij}E_{X^i}
+\partial_i\alpha^{kj}(\eta_k \wedge E_{\eta_i})-dE_{\eta_j}=0.
\label{Noether}
\end{equation}

\section{First steps of the Batalin-Vilkovisky formalism}

Rather than review the Batalin-Vilkovisky formalism in general as in
\cite{HT,bv:closure,bv:closure2}, we illustrate it by example: the
Poisson sigma model we
have been
considering.  Batalin and Vilkovisky first construct a  graded
commutative
algebra
over $Loc_E $ with generators
 $X_i^+$ and $\eta^{+i}$, called `anti-fields', $\gamma_i$ called
`ghosts'
and $\gamma^{+i}$, called `anti-ghosts. (If only the ghosts were used as generators, this would be a BRST algebra.)

These generators are bigraded, as indicated in the following table where the form
degree is displayed as the top row and the ghost degree as the first column.
  {\it The graded commutativity is with respect to the sum of the ghost
degree and
the form degree (which we call the total degree).}

 The
assignments of degree (from left to right) and ghost
number (from top to bottom) are given by
\[
\begin{array}{rccc}
 & 0&1&2                \\
-2&    &    & \bbeta^{+i}    \\
-1&    & \eta^{+i} & X^{+}_i\\
0& X^i   &\eta_i & \\
1& \bbeta_i& &
\end{array}
\]
Ultimately, this algebra is given a  differential $D$ which is a
derivation
with respect to the ghost
degree, but initially has just two such derivations
which need not square to zero.

One of the derivations $\delta$ looks like the Chevalley-Eilenberg
differential for Lie algebra cohomology, even though Batalin
and Vilkovisky need not have a
Lie algebra, and is often called a BRST operator. It is defined
initially by

\begin{align}
\delta X^i&=\alpha^{ij}(X)\bbeta_j,\\
\delta \eta_i&=-d\bbeta_i-
\partial_i\alpha^{jk}(X)\eta_j\bbeta_k, \notag \\
\delta \gamma_i &= \frac12\,
\partial_i\alpha^{jk}(X)\bbeta_j\bbeta_k.\notag
\label{deltafields}
\end{align}
\noindent

The other derivation, $d_{KT},$ does square to zero.  It is
the Koszul-Tate differential for the differential ideal generated
by the Euler-Lagrange equations. The Koszul complex is graded by the
ghost number. This means the
anti-fields generate the Koszul complex with
\def\dkt{d_{KT}}
\def\dce{d_{CE}}
\begin{align}
\dkt \ X^{+}_i&=
d\eta_{i}
+\frac12\,\partial_i\alpha^{kl}(X)\eta_k\wedge\eta_l =  E_{X^i}
\\
\dkt \ \eta^{+i}&=
-dX^i-\alpha^{ij}(X)\eta_j =  E_{\eta_i}.\notag
\label{KTantifields}
\end{align}
Because of the Noether identities,
the Koszul complex has non-trivial cohomology in ghost degree $-1$,
namely the classes given by the formulas for the identities with
$E_{X^i}$ and
$E_{\eta_i}$
replaced by $X^{+}_i$ and $\eta^{+i}:$
\begin{equation}
-\alpha^{ij}(X)X^+_j
-\partial_k\alpha^{ij}(X)\eta_j\wedge \eta^{+k}-d\eta^{+i}.
\end{equation}
These classes can be killed by adjoining the anti-ghosts
$\gamma^{+i}$ and defining
\begin{equation}
\dkt \bbeta^{+i}=
-\alpha^{ij}(X)X^+_j
-\partial_k\alpha^{ij}(X)\eta_j\wedge \eta^{+k}-d\eta^{+i}.
\label{KTantighosts}
\end{equation}
Thus the anti-ghosts occur precisely because of the identities
identified by Noether.

The pairing between symmetries and identities is now expressed as the
pairing between ghosts and anti-ghosts, which plays a crucial role in
the Batalin-Vilkovisky anti-bracket, but first
the anti-fields and anti-ghosts are themselves subject to symmetries
corresponding
to $\delta_\beta$ as follows:

\begin{align}
\delta X_i^+ &= \partial_i\alpha^{kj}(X) X_k^+\gamma_j\\
\delta \eta^{+i} &= \partial_k\alpha^{ij}(X)\eta^{+k}\gamma_j\notag\\
\delta \gamma^{+i} &= \partial_k\alpha^{ij}(X)\gamma^{+k}\gamma_j.\notag
\label{delta-anti}
\end{align}

\section{The Batalin-Vilkovisky anti-bracket and total differential}
The hoped for total differential $D$ will be obtained by adding `terms
of
higher order'
to $\dkt + \delta,$ which does not square to zero.
To do this in general,
Batalin and Vilkovisky introduce an `anti-bracket' $(\ ,\ )$
which is defined in terms of distributional derivatives of functionals
of
the fields and anti-fields.

Before we define the anti-bracket, it is convenient to first consider
the
definition of the derivative of a functional $A$  of fields
and antifields which are denoted collectively as $(\psi_\alpha).$ The
derivative
$\frac{\partial A}{\partial{\psi^\alpha}}$ is the
distribution whose value at test forms  $(\rho^\alpha)$ (of the same
degree and
ghost number as $(\psi^\alpha)$) is given by
\[
\frac d{dt}
A(\psi+t\rho)\biggl|_{t=0}=
\int_{\Sigma}\rho^\alpha\wedge
\frac{\partial A}{\partial{\psi^\alpha}}.\]

Consider  the functional $A$ defined by
$$A(\phi,\phi^+)=\int_{\Sigma} (\phi \wedge \phi^+)$$
then we see that up to signs $\frac {\partial A}{\partial \phi}$ is in
some sense
identified
with $\phi^+$ while $\frac {\partial A}{\partial \phi^+}$ is identified
with
$\phi.$ In this way we see that $\phi$ and $\phi^+$ are ``canonically
conjuguate".

Thus we have a canonical distributional pairing of each field or ghost
with
its `anti':
\begin{align}
(X^i, &X_j^+)= \delta^i_j \\
(\eta_j, &\eta_+^i) \ = \delta^i_j \notag\\
(\bbeta_j, &\bbeta_+^i) \ = \delta^i_j.\notag
\end{align}

The BV anti-bracket extends this as a graded biderivation with repect to
ghost
degree and in this example can be written as $(A,B) =$

\newcommand{\cev}[1]{{\stackrel{\leftarrow}{#1}}}
\begin{eqnarray}
\sum_\alpha
\int_{\Sigma} (-1)^{|\phi_{\alpha}|(|\phi_{\alpha}|+|A|)}\left(
\frac {\partial A}{\partial\phi^\alpha}\wedge
\frac{\partial B}{\partial{\phi^+_\alpha}}
-(-1)^{(deg(\phi_\alpha)+|A|+1)}
\frac {\partial A}{\partial\phi^+_\alpha}\wedge
\frac{\partial B}{\partial{\phi^\alpha}}\right)
\end{eqnarray}
where $|C|=gh(C)+deg(C)$ denotes the Grassman parity of
$C$ ($C$ is either a field or a function of fields). Note that
physicists prefer to
use both left and right derivatives and hence exhibit a different set of
signs.

The antibracket obeys the graded commutativity relation
\[
(A,B)=-(-1)^{(\mathrm{gh}(A)-1)(\mathrm{gh}(B)-1)}(B,A)
\]
and the Leibnitz rule
\begin{equation}\label{e-leibnitz}
(A,BC)=(A,B)C+(-1)^{(\mathrm{gh}(A)-1)\mathrm{gh}(B)}B(A,C),
\end{equation}
which emphasizes the resemblance to a Poisson bracket.  The only
difference from a graded Poisson bracket is that the bracket shifts
the degree by 1
and the several identities (skew-commutativity, Jacobi and
 Leibniz) inherit certain signs. Such an `odd' Poisson bracket is also
known as a {\bf Gerstenhaber bracket} \cite{gerstenhaber:AM63}.

Now it is possible to express $\dkt +\delta$ in the form $(S^0+S^1, \ \
)$
where $$S^0= (X,\eta)=\int_{\Sigma}(\eta_i\wedge dX^i) +
\frac {1}{2} (\alpha^{ij}\circ X)(\eta_i\wedge \eta_j),$$
 our original action, and $S^1$ is
\begin{eqnarray}
\int_{\Sigma}
 X_i^+\alpha^{ij}(X)\gamma_j
-\eta^{+i}\wedge
(d\gamma_i+\partial_i\alpha^{kj}(X)\eta_k\gamma_j)
-\frac12\,\gamma^{+i}
\partial_i\alpha^{jk}(X)\gamma_j\gamma_k .
\end{eqnarray}

Corresponding to the fact that $(\dkt + \delta)^2 \neq 0,$
we have
$$(S^0+S^1, S^0+S^1) \neq 0. $$

The additional terms in the differential $D$ we seek will be found by
extending
$ S^0+S^1$ by terms of higher order to achieve the full BV action
$S_{BV}.$.
First, let us analyze the derivation $(S^0,\ ).$

Notice that $(\ ,X^+_i)$ is effectively (up to sign)
$\partial_{X^i}$ and similarly for the other
anti's, while $(\ ,X^i)$ is effectively $\partial_{X_i^+},$ etc.
More precisely, for any of our basic variables, denoted collectively
as $\phi^a$ and their anti's denoted $\phi^+_a$, we have

\begin{align}
(S ,\phi^a) &= (-1)^{gh (\phi^a)}\frac {\partial S}{\partial \phi^+_a}\\
(S ,\phi^+_a) &= (-1)^{gh (\phi^a)+deg(\phi^a)}\frac {\partial
S}{\partial \phi^a}
\end{align}
whenever the parity of $S(\phi^{\alpha},\phi_{\alpha}^+)$ is even. The
parities
of $S^0,S^1,S^2$ are all 2. Recall that the parity was defined above to
be the total
degree.

Since $S^0$ has no anti's,
$(S^0,S^0)=0,$ in fact,
$$(S^0, X^i) =0,\ \ (S^0,\eta_i) = 0 \ {\rm and}\  (S^0, \gamma_i) =0.$$
\noindent
However, $(S^0,\ )$ does act non-tivially on some of the anti's:
\begin{align}
(S^0, X_i^+) &= d\eta_i + 1/2
\partial_i\alpha^{kl}(X)\eta_k\wedge\eta_j\\
(S^0,\eta^{+i}) &=  -dX^i -\alpha^{ij}\eta_j,\notag
\end{align}
which reproduces part of $\dkt$, cf. (\ref{KTantifields}), while
$(S^0,\gamma^{+i}) = 0.$

Now consider $(S^1,\ ):$
\begin{align}
(S^1, X^i) &= \alpha^{ij}\gamma_j\\
(S^1, \eta_i) &= -(d\gamma_i +
\partial_i\alpha^{jk}(X))\wedge\eta_j\gamma_k\notag\\
(S^1, \gamma_i)&= 1/2 \partial_i\alpha^{jk}(X))\gamma_j\gamma_k\notag\\
(S^1, X^+_i) &= \partial_i\alpha^{kl}(X)X^+_k\gamma_l
-\partial_i\partial_j\alpha^{kl}(X)
\eta^{+j}\wedge\eta_k\gamma_l
-\frac12\,\partial_i\partial_j\alpha^{kl}(X)\gamma^{+j}\gamma_k\gamma_l
\notag\\
(S^1,\eta^{+i}) &= \eta^{+k}\partial_k\alpha^{ij}(X)\gamma_j\notag\\
(S^1,\gamma^{+i}) &=-\alpha^{ij}X_i^+ -d\eta^{+i}
 +\partial_k\alpha^{ij}(X)\eta^{+k}\wedge\eta_j
+\partial_k\alpha^{ij}(X)\gamma^{+k}\gamma_j,\notag
\end{align}
reproducing (\ref{deltafields}) and (\ref{delta-anti}).

Batalin and Vilkovisky show that, in much more general situations,
one can add terms $S^i$ of ghost degree $i>1$  to achieve a total
$S_{BV}$
such that
$$(S_{BV}, S_{BV}) = 0.$$
The reason for this is that the $\dkt$ homology vanishes in appropriate
degrees.

In the Cattaneo-Felder model, only one more term is needed:
\begin{eqnarray}
S^2 = \int_{\Sigma} -\frac14\,\eta^{+i}\wedge\eta^{+j}
\partial_i\partial_j\alpha^{kl}(X)\gamma_k\gamma_l.
\end{eqnarray}
Thus the total Batalin-Vilkovisky generator is
\begin{eqnarray}
S_{BV}
&=&
\int_{\Sigma}
\eta_i\wedge dX^i
+\frac12\alpha^{ij}(X)\eta_i\wedge\eta_j\\
&&
+ X_i^+\alpha^{ij}(X)\gamma_j
-\eta^{+i}\wedge(d\gamma_i+\partial_i\alpha^{kl}(X)\eta_k\gamma_l)
-\frac12\,\gamma^{i}_+
\partial_i\alpha^{jk}(X)\gamma_j\gamma_k \notag\\
&&
-\frac14\,\eta^{+i}\wedge\eta^{+j}
\partial_i\partial_j\alpha^{kl}(X)\gamma_k\gamma_l.\notag
\end{eqnarray}

\section{Summary}
We hope to have called deserving attention to Noether's second
variational
theorem and how it accounts for the anti-ghosts which are an essential
part
of the Batalin-Vilkovisky method. Beyond that, we are now able to show
how the terms $S^i$ in the total $S_{BV}$
of the Catanneo-Felder sigma model correspond to the Koszul-Tate,
Chevalley-Eilenberg  and other parts of the total differential
in the BV differential graded algebra.  Consider the total differential
as found in Cattaneo and Felder:
\begin{eqnarray}
\delta X^i&=&\alpha^{ij}(X)\bbeta_j,
\\
\delta\eta^{+i}&=&
-dX^i-\alpha^{ij}(X)\eta_j+
\partial_k\alpha^{ij}(X)\eta^{+k}\bbeta_j,\\
\delta\bbeta^{+i}&=&
-d\eta^{+i}
-\alpha^{ij}(X)X^+_j
+\frac12\,\partial_k\partial_l\alpha^{ij}(X)
\eta^{+k}\wedge\eta^{+l}\bbeta_j
\\
 &&+\partial_k\alpha^{ij}(X)\eta^{+k}\wedge\eta_j
+\partial_k\alpha^{ij}(X)\bbeta^{+k}\bbeta_j.
\end{eqnarray}
and
\begin{eqnarray}
\delta \bbeta_i&=&\frac12\,\partial_i\alpha^{kl}(X)\bbeta_k
\bbeta_l,
\\
\delta\eta_i&=&
-d\bbeta_i-\partial_i\alpha^{kl}(X)\eta_k\bbeta_l-
\frac12\,\partial_i\partial_j\alpha^{kl}(X)\eta^{+j}\bbeta_k
\bbeta_l,\\
\delta X^{+}_i&=&
d\eta_{i}
+\partial_i\alpha^{kl}(X)X^+_k\bbeta_l
-\partial_i\partial_j\alpha^{kl}(X)
\eta^{+j}\wedge\eta_k\bbeta_l
+\frac12\,\partial_i\alpha^{kl}(X)\eta_k\wedge\eta_l
\notag\\
 &&-\frac14\,\partial_i\partial_j\partial_p
\alpha^{kl}(X)\eta^{+j}\wedge\eta^{+p}\bbeta_k\bbeta_l
-\frac12\,\partial_i\partial_j\alpha^{kl}(X)\bbeta^{+j}\bbeta_k
\bbeta_l.
\end{eqnarray}

These individual terms can be identified as coming from a particular
$S^i.$ For example,
$\delta X^i$ comes from $S^1$, the first two terms of
$\delta\eta^{+i}$ come from $S^0$ and the third from $S^1$,
as do all the terms of $\delta\bbeta^{+i}$ except for the middle term
which comes from $S^2.$
Similarly,
$\delta \bbeta_i$ comes from $S^1$, the first two terms of
$\delta\eta_i$ come from $S^1$ and the third from $S^2$,
while the five terms of $\delta X^{+}_i$ come from $S^i$ with
$i$ respectively $0,1,1,0,2,1.$

In contrast, if we identify terms as coming from $\dkt$ or $\dce$
we find
$\delta X^i$ comes from $\dce$, the first two terms of
$\delta\eta^{+i}$ come from $\dkt$ and the third from $\dce$, while
the first, second and fourth terms of $\delta\bbeta^{+i}$ come from
$\dkt$,
the fifth from $\dce$ and the third term is of neither origin.
Similarly,
$\delta \bbeta_i$ comes from $\dce$, as do the first two terms of
$\delta\eta_i$   and the third
is of neither origin. The first and fourth terms of $\delta X^{+}_i$
come from
$\dkt,$ the second term comes  $\dce,$ and the remaining terms come from
neither
$ \dce$ nor $\dkt.$

\end{document}